\documentclass[%
reprint,
superscriptaddress,
showpacs,
hidelinks,
amsmath,amssymb,
aps,
prb,
]{revtex4-1}
\usepackage{graphicx,epstopdf}
\usepackage{placeins}
\usepackage{mathtools}
\usepackage{tikz}
\usepackage{tikz-3dplot}
\usetikzlibrary{calc,math}
\usetikzlibrary{intersections}
\usetikzlibrary{positioning}
\usetikzlibrary{arrows,decorations.markings}
\usepackage{pgfplots}
\pgfplotsset{compat=newest}
\usepgflibrary{shapes.geometric}
\usepackage{mathrsfs}
\usepackage{bm}
\usepackage{esint}
\usepackage{siunitx}
\usepackage{physics}
\usepackage{circledsteps}
\usetikzlibrary{arrows.meta}
\usepackage{relsize}

\begin{document}
\title{Visualizing Stokes' theorem with Geometric Algebra}
\author{Kristjan Ottar Klausen}

\address{School of Technology, Reykjavik University, Menntavegur 1, IS-101 Reykjavik, Iceland.}

\vspace{10 mm}

\begin{abstract}
Foundational cases of the generalized Stokes' theorem are visualized using geometric algebra. From considering bivector valued fields, two seldom used instances of the theorem are obtained. Graphical representations are given, showing a radial spin structure corresponding to the generalized curl and a toroidal vortex for the generalized divergence.
Dualities and applications to natural phenomena with non-trivial topology are discussed.
\end{abstract}
\vspace{-20 mm}
\maketitle
\section{INTRODUCTION}
In the past century, Clifford algebras have steadily gained popularity and  found many applications across diverse fields of science in recent years \cite{Hitzer2022,Eduardo2021}. An attractive feature of Clifford algebras is that they unify and generalize various branches of mathematics commonly applied in physics \cite{Hestenes_Stokes_split,Baylis1996}. 
Importantly, Clifford algebras lie at the heart of the recently formulated periodic table of topological invariants\cite{Periodic_Kitaev, Periodic_Classification_Shinsei,Periodic_Sobczyk}, applicable to topological insulators and topological superconductors \cite{Sato_2017,Qi_2011,Hasan_2010}.
Systems with topological phases of matter, offer the possibility of topological invariants in the energy spectrum and play a prominent part in the current race for quantum computation \cite{nature_news}. By closing and reopening of the energy spectra, edge states emerge such as Majorana zero modes \cite{SchapersGroup2019, Review_Agauado,Review_Beenakker} which can function as qubits, robust to environmental perturbations \cite{MajoranaQubit,Kitaev}. Furthermore, within Clifford algebras, general operations of quantum computing and logic can be worked with \cite{Cafaro_2011,HavelDoran2000,quantumerr_cliff_2000,Somaroo_1998}.\\

Geometric algebra, a real valued formulation of Clifford algebra, has proven to be a straight-forward and intuitive generalization of vector algebra with definite geometric interpretations \cite{doran_lasenby_2003} and clear relations to linear algebra \cite{MacDonald2010_linear,Linear_Sobczyk}.
Geometric algebra can be further formulated to geometric calculus, in which the \textit{boundary theorem} can be proven, which generalizes the fundamental theorems of scalar and vector calculus to $n-$dimensions \cite{MacDonald2012} for $n$ in $\mathbb{N}$.  The boundary theorem encompasses the generalized Stokes theorem written with differential forms along with its dual cases \cite{Hestenes_DF}, enables the definition of coordinate free integration \cite{Alho2017} and from it, an $n-$dimensional generalization of the residue theorem can be obtained \cite{Hestenes_Sobczyk}.
In light of the above considerations, the value of standardizing geometric algebra in the current mathematics and physics curriculum is apparent.
 Geometric algebra in no way invalidates other fields of mathematics, on the contrary, it unifies and brings fresh perspectives to multiple mathematical disciplines.  Many properties of the algebra have been derived within other formalisms \cite {Hestenes_Sobczyk}. The main advantage of geometric algebra is accessibility and ease of visualization.\\
 
In this work, significant cases of the boundary theorem are visualized using geometric algebra and interrelations to standard theorems of scalar and vector calculus, and their generalizations, clarified. Manifestations in the natural world are discussed.
Solely the necessary components for considering the boundary theorem are introduced as multiple textbooks exist covering the the basics of geometric algebra and calculus \cite{doran_lasenby_2003,MacDonald2010_linear,MacDonald2012,HestenesNFCM,Hestenes_Sobczyk}.

\section{Stokes' theorem}
The fundamental theorem of calculus relates the derivative of a function to its integral and is the crux of calculus,
\begin{equation}
\int_{a}^{b} f'(x) \, dx = f(b) - f(a),
\label{fundam}
\end{equation}
where $f(x)$ is continuously differentiable real-valued function on the open interval (a,b).	
The theorem allows for 
calculation of the integral of a function by evaluating the function's antiderivative at the end points, given that the antiderivative is expressible in terms of elementary functions.

In vector calculus, the theorem has generalizations in two- and three-dimensional spaces using fields along with the vector differential operator $\boldsymbol{\nabla}$.
A vector field is said to be conservative if it can be written as the gradient of a function,
$\mathbf{E}=\boldsymbol{\nabla} \phi $, where $\phi$ is termed the scalar potential for $\mathbf{E}$.
The gradient theorem states that an integral of a conservative field along a path equals the difference in the value of the potential field at the end points,
\begin{equation}
\int_{a}^{b} \boldsymbol{\nabla} \phi \cdot d\boldsymbol{\ell} = \phi(b) - \phi(a).
\label{grad}
\end{equation}
Integration along a closed path $C$ in this case equals zero,
\begin{equation}
	\oint_C \boldsymbol{\nabla} \phi = 0.
\end{equation}
Furthermore, every conservative field is irrotational, as the curl vanishes according to 
\begin{equation}
	\boldsymbol{\nabla} \times \mathbf{E}=\boldsymbol{\nabla} \times (\boldsymbol{\nabla} \phi)=0. 
\end{equation}
The converse is only true if there are no holes in the space in which the integral over $C$ is evaluated, in formal terms if $C$ is simply connected or more generally, locally contractible \cite{Burns2019, Nakahara}.  By considering subsets of spaces with parts removed, non-conservative irrotational fields can be formulated. This study further leads to the mathematics of de Rham cohomology \cite{Nakahara}.\\

A vector field is said to be solenoidal if it can be written as the curl of another vector field, $\mathbf{B} = \boldsymbol{\nabla} \times \mathbf{A}$, where $\mathbf{A}$ is termed the vector potential for $\mathbf{B}$.
In the same way that the integral of the differentiated potential in Eq. \eqref{grad} was given by its value at the boundary, Green's curl theorem states that integration along a closed path forming the boundary of a surface, $C=\partial S$, is non-zero and given by the total curl (circulation density) on the surface,
\begin{equation}
	\iint_S \boldsymbol{\nabla} \times \mathbf{A} \cdot d\boldsymbol{S} = \oint_{\partial S} \mathbf{A} \cdot d\boldsymbol{\ell}.
	\label{curl}
\end{equation}
Furthermore, every solenoidal field is incompressible, without sources/sinks, as the divergence vanishes according to
\begin{equation}
	\boldsymbol{\nabla} \cdot \mathbf{B} = \boldsymbol{\nabla} \cdot (\boldsymbol{\nabla} \times \mathbf{A})=0.
\end{equation}
Any smooth vector field $\mathbf{F}$ (vanishing at infinity) in $\mathbb{R}^3$ can be written as the sum of an irrotational and a solenoidal field via Helmholtz decomposition \cite{HH_survey,HH_OV,HH_MM}
\begin{equation}
	\mathbf{F}= \mathbf{E}+\mathbf{B} = \boldsymbol{\nabla} \phi + \boldsymbol{\nabla} \times \mathbf{A} .
	\label{HH}
\end{equation}
Although the divergence of the solenoidal component vanishes, the divergence of the irrotational component does not, and is given by the Laplacian of the scalar potential,
\begin{equation}
	\boldsymbol{\nabla} \cdot \mathbf{E}= \boldsymbol{\nabla} \cdot (\boldsymbol{\nabla} \phi) = \boldsymbol{\nabla}^2 \phi .
\end{equation}
The divergence theorem equates the volume integration of this expression to the flux of the field $\mathbf{E}$ through a surface $S=\partial V,$ containing the sources/sinks of $\mathbf{E}$,
\begin{equation}
\iiint_V \boldsymbol{\nabla} \cdot \mathbf{E} \, dV = \oiint_{\partial V} \mathbf{E} \cdot d\boldsymbol{S}.
\label{div}
\end{equation}
Utilization of the above theorems is central to electromagnetism and therefore notation has been chosen such that results from electro- and magnetostatics can be recognized. 
The theorems from Eqs.\ \eqref{grad},\eqref{curl},\eqref{div} can be combined and generalized to $n$-dimensions, $n \in \mathbb{N}$, in a single theorem known as Stokes' theorem using differential forms along with the exterior derivative \cite{Hestenes_Stokes_split,Hestenes_DF, Cartan}. Although attributed to G.G.\ Stokes, the theorem is the result of a cumulative effort of multiple researches in the past two centuries \cite{Katz79}. Before discussing Stokes' theorem in geometric calculus, let us briefly introduce the basics and structure of geometric algebra.

\section{Geometric algebra}
From the concept of vectors, geometric algebra can be built up using the geometric product, which defines vector multiplication and is given by
\begin{equation}
	\mathbf{ab= a\cdot b + a\wedge b} \ .
	\label{geoprod}
\end{equation}
The center dot denotes a symmetric inner (scalar) product,

\begin{equation}
	\mathbf{a\cdot b} = \frac{1}{2} (\mathbf{ab} + \mathbf{ba}).
\end{equation}
The wedge operator $\wedge$,
denotes the anti-symmetric outer (wedge) product,
\begin{equation}
\mathbf{a} \wedge \mathbf{b}= \frac{1}{2} (\mathbf{ab} - \mathbf{ba}).
\end{equation}

\begin{figure}[h!]
	\centering
	\includegraphics[scale=1]{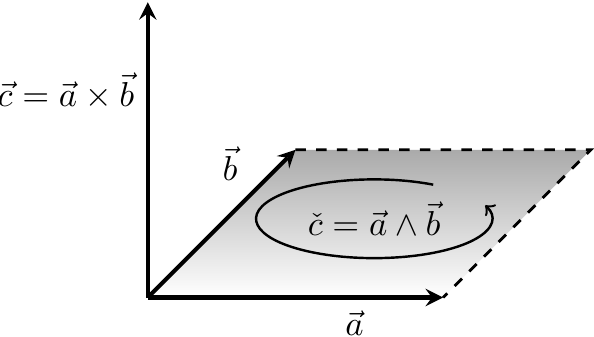}
	\caption{The vectors $\vec{a}$ and $\vec{b}$ span the plane of the bivector $\check{c}$. The vector $\vec{c}$ is orthogonal to the plane and in three-dimensional space only, dual to the bivector.}
	\label{crosspwedge}
\end{figure}

The outer product is a generalization of the cross product to $n$-dimensions yet still not equal to the cross product in three dimensions, but dual to it. For two orthogonal vectors $\vec{a}\text{ and }\vec{b}$, the cross product $\vec{a} \times \vec{b} = \vec{c}$ is a third vector, with magnitude and direction, orthogonal to both $\vec{a}\text{ and }\vec{b}$. The wedge product $ \vec{a} \wedge \vec{b} = \check{c}$ is however a two-dimensional bivector lying in the plane spanned by $\vec{a}\text{ and }\vec{b}$, having magnitude and area, Fig.\ \ref{crosspwedge}.
In this way the outer product is dimension-raising, whilst the inner product is dimension-reducing. In general the outer product between an $n$-vector and $m$-vector is an $p$-vector of dimension $p=n+m$ whilst their
	inner product is an $l$-vector of dimension $l=n-m$ for $n\ge m$. It's important to note that the inner product can be defined in a more detailed manner with the concept of contractions \cite{Dorst2002}.

Using the outer product, higher-dimensional (higher grade) vectors can be defined such as three-dimensional trivectors $ \tilde{d} = \vec{a} \wedge \vec{b} \wedge \vec{c} ,$ four-dimensional quadvectors and in general $n$-vectors (termed blades in Clifford algebra) which form the building blocks of geometric algebra. 
Blades of different grades can be added resulting in linear combinations termed multivectors, which can fairly be said to generalize the concept of number to higher dimensions \cite{Periodic_Sobczyk}.

Any $n$-vector $\mathbf{v}$ determines an $n$-dimensional vector space $\mathcal{V}_n$ of all vectors $\mathbf{x}$ satisfying $\mathbf{x}\wedge \mathbf{v}=0$. Addition and multiplication of vectors in $\mathcal{V}_n$ results in a $2^n$-dimensional linear space closed under the geometric product, the geometric algebra $\mathbb{G}_n$. 
For orthogonal base vectors $\mathbf{e}_1$ and $\mathbf{e}_2$ we have the convenient notation for their geometric product
	\begin{equation}
		\mathbf{e}_1 \mathbf{e}_2 = \mathbf{e}_1 \cdot \mathbf{e}_2 + \mathbf{e}_1 \wedge \mathbf{e}_2 = 0 + \mathbf{e}_1 \wedge \mathbf{e}_2 = \mathbf{e}_{12}.
\end{equation}
Since an $n$-vector in $n$-dimensional space does not have degrees of freedom for orientation, other that up to a sign, such an element becomes a \textit{pseudoscalar} \cite{Hestenes_Sobczyk}.  As an example, the complex number $i$ is a unit bivector which becomes a pseudoscalar in the 2D plane. In this way, complex numbers can be understood as the even part of $\mathbb{G}_2$ geometric algebra\cite{imag}.

The binomial expansion of geometric algebra $\mathbb{G}_n$ with unit vectors $e_n^2=1$ for $n=0,1,2,3,4$ is shown in Fig.\ \ref{binom} along with the relations of its even components to the division algebras of complex numbers, quaternions and octonions $\mathbb{C}, \mathbb{H} \text{ and } \mathbb{O}$, respectively.
$\mathbb{G}_0$ is the algebra of the real numbers $\mathbb{R}$, the only other division algebra. In $\mathbb{G}_1$, vectors are pseudoscalars and so the geometric algebra is essentially the algebra of the real number line, $\mathbb{R}^1$.  For the plane,
$\mathbb{G}_2$ consists of scalars, two orthogonal vectors and a pseudoscalar bivector. $\mathbb{G}_3$ consists of scalars, three orthogonal vectors, three orthogonal bivectors and a pseudoscalar trivector. $\mathbb{G}_4$ consists of scalars, four orthogonal vectors, six orthogonal bivectors, four trivectors and a pseudoscalar quadvector.
As the product of even multivectors is always even\cite{Hestenes_Sobczyk}, the even components of geometric algebra in each dimension form a set closed under multiplication, a subalgebra of rotations. The division algebras have been shown to form a basis for all compact symmetric spaces \cite{Huang2011}.\\

\begin{figure}[t]
	\centering
	\begin{tabular}{rccccccccc}
		$\mathbb{G}_0$:&    &    &    &    &  1\\\noalign{\smallskip\smallskip}
		$\mathbb{G}_1$:&    &    &    &  1 &    &  1\\\noalign{\smallskip\smallskip}
		$\mathbb{C} \in \mathbb{G}_2$:&    &    & \Circled{1} &    &  2 &    &  \Circled{1} \\\noalign{\smallskip\smallskip}
		$\mathbb{H} \in \mathbb{G}_3$:&    &   \Circled{1}&    &  3 &    &   \Circled{3} &    &  1\\\noalign{\smallskip\smallskip}
		$\mathbb{O} \in \mathbb{G}_4$:&  \Circled{1}&    &  4 &    &  \Circled{6} &    &  4 &    &  \Circled{1} \\\noalign{\smallskip\smallskip}
	\end{tabular}
	\caption{Binomial expansion of geometric algebra up to $\mathbb{G}_4$, showing the hypercomplex division algebras as the even components (circled) in each dimension.}
	\label{binom}
\end{figure}
The geometric product in Eq.\ \eqref{geoprod} is valid for the product of a vector and multivector\cite{MacDonald2012}, but the geometric product of two multivectors generally involves more terms \cite{Hestenes_Sobczyk}. The product of two bivectors $A$ and $B$ has three terms: scalar, bivector and quadvector,
\begin{equation}
AB=A \cdot B + A \times B + A \wedge B.
\end{equation}
The quadvector term is non-zero in $\mathbb{G}_4$ and above. The bivector component is denoted with the grade-preserving antisymmetric commutator product 
\begin{equation}
A \times B = \frac{1}{2} [AB -BA].
\label{com}
\end{equation}
Using the same notation for the commutator product and the cross product could be seen as a source of confusion, but keeping in mind the grade of the entities under multiplication clears this up as the cross product is only defined for vectors. The similar properties of the two products is due to the duality of vectors and their orthogonal planes in three dimensional space.

\section{Geometric calculus and the boundary theorem}
The multivector derivative operator $\boldsymbol{\nabla}$ for a differentiable multivector field $F$ with an orthonormal basis $\mathbf{e}_i$ on an open set in $\mathbb{R}^n$ can be defined by the following expression, using the shorthand notation $\partial_i = \frac{\partial}{\partial_i}$,
\begin{equation}
\boldsymbol{\nabla} F = \mathbf{e}_i \partial_i F = \mathbf{e}_1 \partial_1 F + \mathbf{e}_2 \partial_2 F + ... + \mathbf{e}_n \partial_n F.
\end{equation}
Considering that the operator $\boldsymbol{\nabla} = \mathbf{e}_i \partial_i$ functions as a vector algebraically, for homogeneous multivector fields it can be factored into two components using the geometric product from Eq.\ \eqref{geoprod}, 
\begin{equation}
\boldsymbol{\nabla} F = \boldsymbol{\nabla} \cdot F + \boldsymbol{\nabla} \wedge F,
\label{nablaF}
\end{equation}
where $\boldsymbol{\nabla} \cdot F$ is the generalized divergence and $\boldsymbol{\nabla} \wedge F$ is the generalized curl. 
In order to differentiate fields on manifolds in $\mathbb{R}^n$, the component of the vector derivative lying within the manifold has to be used since in general it has components parallel and perpendicular to the manifold\cite{MacDonald2012}, 
\begin{equation}
	\boldsymbol{\nabla} = \boldsymbol{\nabla}_{||}  + \boldsymbol{\nabla}_\perp,
	\label{vd_proj}
\end{equation} respectively.
Let $M^m$ be a smooth $m$-dimensional oriented manifold in $\mathbb{R}^n$ with a piece-wise smooth boundary $\partial M$ of integer dimension $m-1$ where $m \ge 1$. The component of the vector derivative that is parallel (intrinsic) to the manifold is given by
\begin{equation}
	\boldsymbol{\nabla}_{||}F= I_m^{-1} (I_m \cdot \boldsymbol{\nabla} F),
	\label{vd_tan}
\end{equation}
where $I_m$ denotes the pseudoscalar of the manifold $M$. Multiplication by $I_m^{-1}$ gives the dual, which is necessary to preserve the grade and directionality. This constitutes the definition of the tangential vector derivative on manifolds 
\begin{equation}
	\boldsymbol{\partial} = \boldsymbol{\nabla}_{||} =  I_m^{-1} (I_m \cdot \boldsymbol{\nabla} ),
	\label{vd_par}
\end{equation}
where the boldface vector notation distinguishes it from the boundary operator $\partial$. For the special case of $m=n$, i.e.\ a manifold $M^n$ in $\mathbb{R}^n$, the vector derivative lies entirely within the manifold so that $\boldsymbol{\nabla}_\perp=0$ and thus $\boldsymbol{\partial}=\boldsymbol{\nabla}$, from Eqs.\ \eqref{vd_proj} and \eqref{vd_par}.

 By defining a directed element of the manifold $M$ as a differential $m$-vector d$\mathbf{x}^m$, and similarly d$\mathbf{x}^{m-1}$ for $\partial M$, the boundary theorem of geometric calculus for flat manifolds\cite{MacDonald2012,Hestenes_STA,Hestenes_Sobczyk} can be written as
\begin{equation}
	\int_{M} \text{d}\mathbf{x}^m \, \boldsymbol{\partial} F = \oint_{\partial M} d\mathbf{x}^{m-1} \, F \ .
\end{equation}
Intuitively, this means that the total change in the field $F$ within a space $M$, equals the state or value of $F$ at the boundary of the space. In other words, information about the change in the field within the space $M$ is encoded on its boundary. Hestenes\cite{Hestenes_DF} noted that from the boundary theorem, the tangential vector derivative can be given a coordinate-free definition,
	\begin{equation}
		\boldsymbol{\partial}F = \lim_{|\text{d}\mathbf{x}^m| \to 0} \frac{I_m^{-1}}{|\text{d}\mathbf{x}^m|} \oint_{\partial M} \text{d}\mathbf{x}^{m-1} F.
		\label{tangvecd}
	\end{equation}

Using the orthogonal decomposition from Eq.\ \eqref{nablaF}, the boundary theorem can be seen to include separate components of generalized divergence and curl\cite{Hestenes_Stokes_split}
\begin{equation}
	\int_{M} \text{d}\mathbf{x}^m \, (\boldsymbol{\partial} \cdot F + \boldsymbol{\partial} \wedge F) = \oint_{\partial M} d\mathbf{x}^{m-1} \, F \ .
	\label{bt_comp}
\end{equation}
The curl component can be considered analogous to the generalized Stokes' theorem for differential forms\cite{Hestenes_DF},
\begin{equation}
	\int_M d\omega = \oint_{\partial M} \omega,
	\label{st_diff}
\end{equation}
 with $\omega=d\mathbf{x}^{m-1}\cdot F$ and $d\omega= \text{d}\mathbf{x}^m \cdot (\boldsymbol{\partial} \wedge F)$. From here on out it will be referred to as the outer part of Eq.\ \eqref{bt_comp}, due to the outer product of the vector derivative with the field $F$. The divergence component is its orthogonal counterpart which completes the boundary theorem and will be referred to as the inner part\cite{Hestenes_Stokes_split}.
The generalized Stokes' theorem of forms, Eq.\ \eqref{st_diff}, is central to cohomology and multiple other aspects of algebraic topology \cite{bott2013differential}. The theorem holds not only for smooth manifolds, but also for $k-$chains (triangulations of $M$) and even non-smooth chains (fractals) \cite{Harrison93}.
From the perspective of geometric algebra however, Stokes's theorem of forms only contains half of the information content of the boundary theorem in geometric calculus \cite{Hestenes_Stokes_split}. The other half is the dual, which is implicit in the geometric product.\\
\section{VISUALIZING THE BOUNDARY THEOREM}
 In order to cover the main cases of interest, the instances for $n=1,2,3$ will be explored where the multivector field $F$ is of single grade and $\dim(F)=\dim(M)-1=dim(\partial M)$.
  The graphical representation is defined in a standard manner, denoting scalars with dots, vectors with oriented lines, bivectors with oriented planes, trivectors with oriented volumes and so forth. For simplicity and familiarity, we take the manifold $M^n$ to be in $\mathbb{R}^n$ unless stated otherwise, so that $\boldsymbol{\partial}=\boldsymbol{\nabla}$. Furthermore, manifolds are chosen to be of maximal symmetry (spherical). The theorem is nonetheless equally applicable for irregular manifolds, as well as manifolds with holes, using a similar treatment as standardized in vector calculus \cite{MacDonald2012}. 
  
\begin{figure}[h!]
 	\includegraphics[width=6cm]{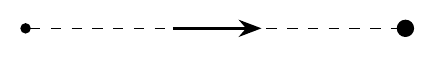}
 	\caption{Graphical representation of the outer part of the boundary theorem for the case of a scalar field on an open one-dimensional manifold (line) $M^1$. The scalar values of $F$ at the boundary $\partial M$ are denoted by points of unequal size, which corresponds to the gradient $\boldsymbol{\nabla}F=\boldsymbol{\nabla} \wedge F$ on $M$, represented by an arrow. The dashed line shows the direction of projection toward the boundary, and coincides with $M^1$.}
 	\label{fund}
 \end{figure}

 In Fig.\ \ref{fund}, a graphical representation of the outer part of the boundary theorem on a one-dimensional manifold is presented, which is essentially the fundamental theorem of calculus, Eq.\ \eqref{fundam} or the gradient theorem for a simple path\cite{MacDonald2012}, Eq.\ \eqref{grad}. 
\FloatBarrier
The inner part in this case is zero
since if the field $F$ is scalar valued then $\boldsymbol{\nabla} \cdot F=0$. Adhering to the dimension reducing property of the inner product results in $dim(\boldsymbol{\nabla} \cdot F)=-1$. Negative-dimensional fields are undefined but may have little explored relations to fractal dimensions \cite{Mandelbrot09,Mandelbrot91}.
 From a physical point of view, the gradient theorem can be understood to state that a difference in scalar values at the boundary of a line will result in a unidirectional flow between them, and vice versa. For example, fluid and charge flow from velocity and electric potential gradients, respectively. The scalar gradient is denoted with points of unequal magnitude at the boundary. The more classical graphical representation is denoting the scalar magnitude difference by a slope, but this simplified representation carries nicely on to higher-dimensions.

For a surface $M^2$ and vector field $\boldsymbol{F}$, the boundary theorem can be split into two orthogonal (dual) components,
\begin{equation}
	\iint_{M}\text{d}\mathbf{x}^2 \cdot (\boldsymbol{\nabla} \wedge \mathbf{F})= \oint_{\partial M} d\mathbf{x}^{1}\cdot \mathbf{F},
	\label{outer1}
\end{equation}
and
\begin{equation}
	\iint_{M}\text{d}\mathbf{x}^2 \wedge (\boldsymbol{\nabla} \cdot \mathbf{F}   ) = \oint_{\partial M} d\mathbf{x}^{1}\wedge \mathbf{F} .
	\label{inner1}
\end{equation}
The outer part, Eq.\ \eqref{outer1}, takes on the familiar form of Green's curl theorem, Eq.\ \eqref{curl}, using that the curl is dual to the outer product
\begin{equation}
	\boldsymbol{\nabla} \wedge \mathbf{F}=I \,\boldsymbol{\nabla} \times \mathbf{F},
\end{equation}
where $I=\mathbf{e}_1 \wedge \mathbf{e}_2 \wedge \mathbf{e}_3$ is the $\mathbb{G}_3$ pseudoscalar, along with $d\mathbf{x}^2= i |dS|= I \mathbf{e}_3 |dS|= I \hat{\mathbf{n}} |dS|$. 
\begin{figure}[h!]
	\includegraphics[height=3.9cm]{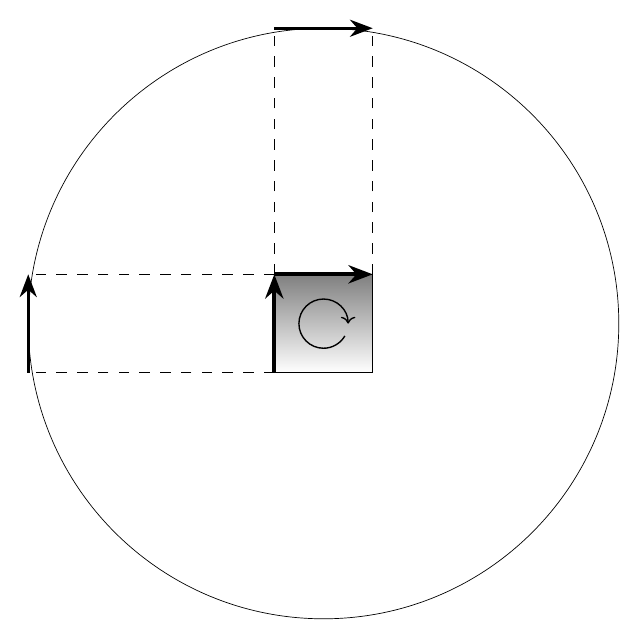}
	\caption{Outer part of the boundary theorem for the surface $M^2$. A non-zero rotation of the vector field $\boldsymbol{F}$, denoted with arrows along the one-dimensional boundary $\partial M$, results from the bivector $\boldsymbol{\nabla} \wedge \mathbf{F}$ within $M$, and vice versa. Orientation is determined by the sign.}
	\label{2Dgreen}
\end{figure}
A graphical representation can be motivated by considering the curl component of Eq.\ \eqref{bt_comp}.
In essence, a bivector leads to a rotation along a closed boundary containing it, measured by the component of $\mathbf{F}$ parallel with the line element $d\mathbf{x}^1$, Fig.\ \ref{2Dgreen}. Conversely, a rotation along a closed curve in the plane is due to a bivector enclosed within it. It is common to convey the idea behind Stokes' theorem by showing how multiple closed loop integral contours cancel within a surface, leaving the boundary term \cite{MacDonald2012}.  In Fig.\ \ref{2Dgreen} however, the bivector $\boldsymbol{\nabla}\wedge \boldsymbol{F}$ is centered to emphasize it's nature as a circulation singularity, as in the physical case of an irrotational vortex\cite{BerryVortex}, and the mathematical case of complex residues \cite{doran_lasenby_2003}.

\begin{figure}[h!]
	\includegraphics[height=4cm]{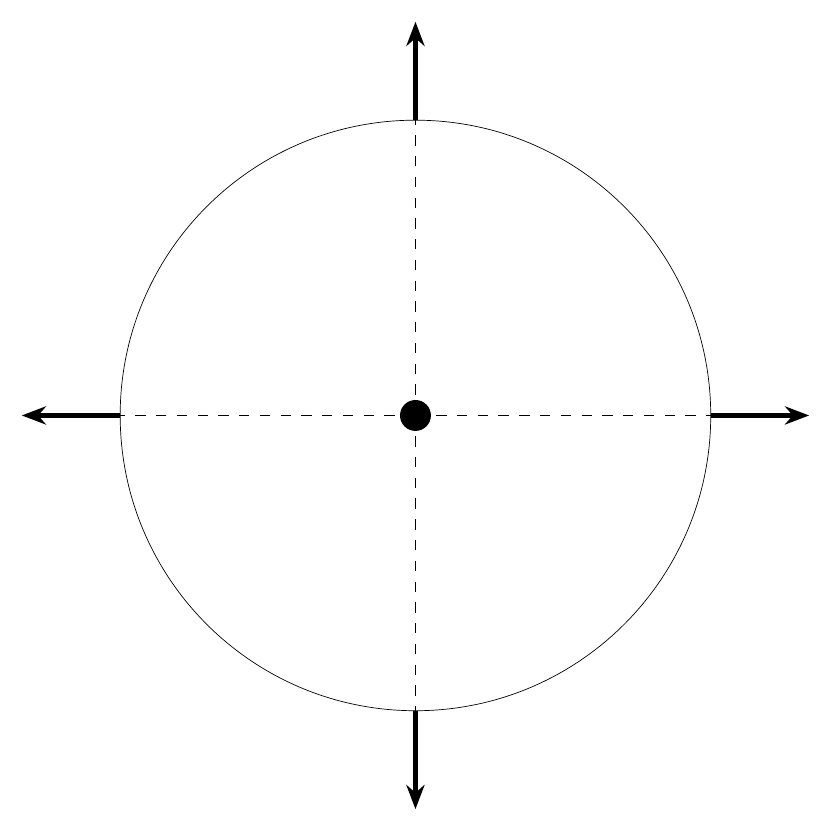}
	\caption{The divergence theorem in two dimensions, the inner part of the boundary theorem for the surface $M^2$.  A non-zero flux, outward (inward), perpendicular to the one-dimensional boundary $\partial M$ stems from the scalar divergence (convergence).}
	\label{2Ddiv}
\end{figure}

Pleasingly, the inner part of the boundary theorem for the surface $M^2$ and vector field $\mathbf{F}$, Eq.\ \eqref{inner1}, describes a two-dimensional divergence theorem, Fig.\ \ref{2Ddiv}. Normally, the divergence is highlighted with more radial rays graphically, but takes on the suggested form here, by adhering to the lowest number of elements needed to span $M$.

\FloatBarrier

\begin{figure}[h!]
 \includegraphics[width=5.5cm]{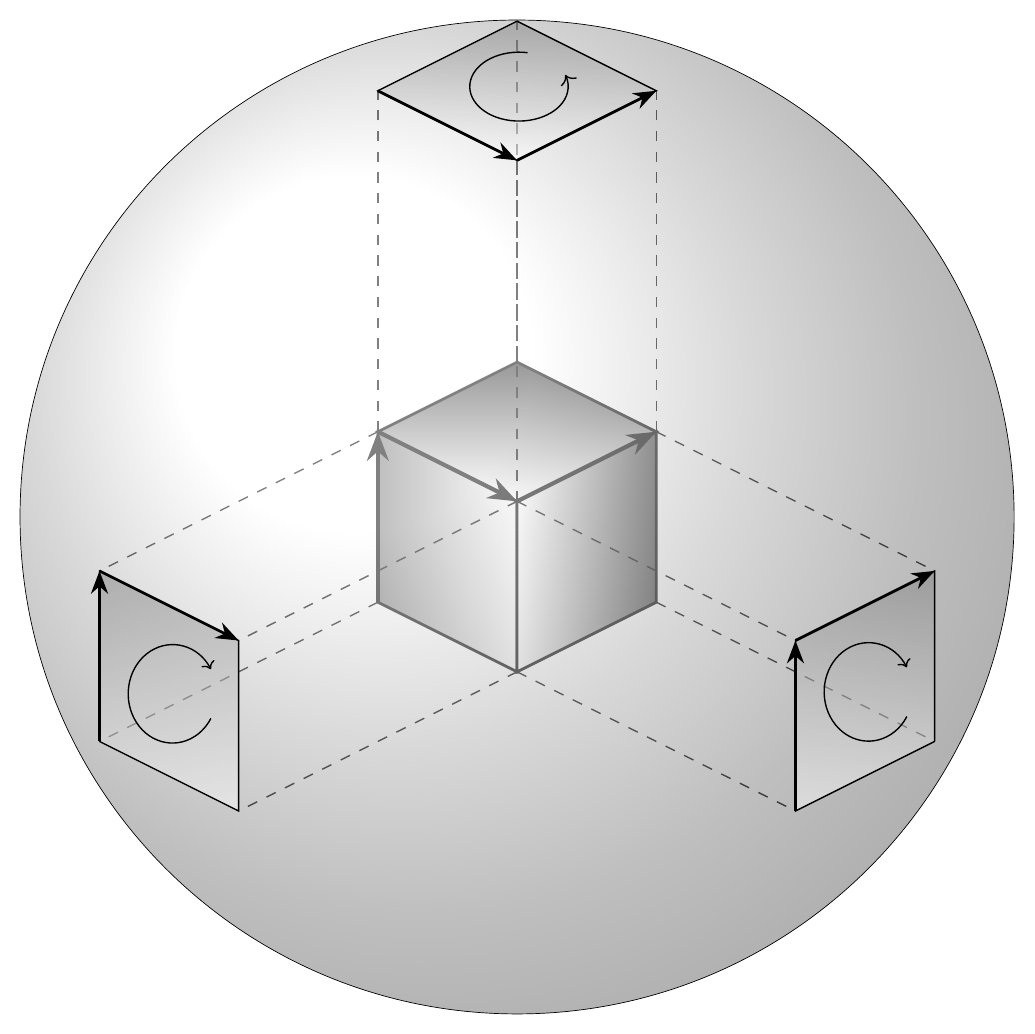}
	\caption{Outer part of the  boundary theorem for the volume $M^3$. The pseudoscalar number density of the trivector valued $\boldsymbol{\nabla} \wedge \check{F}$ equals the sum of the tangent bivector field on the boundary.}
	\label{3Douter}
\end{figure}

\begin{figure}[h!]
	\includegraphics[width=5.65cm]{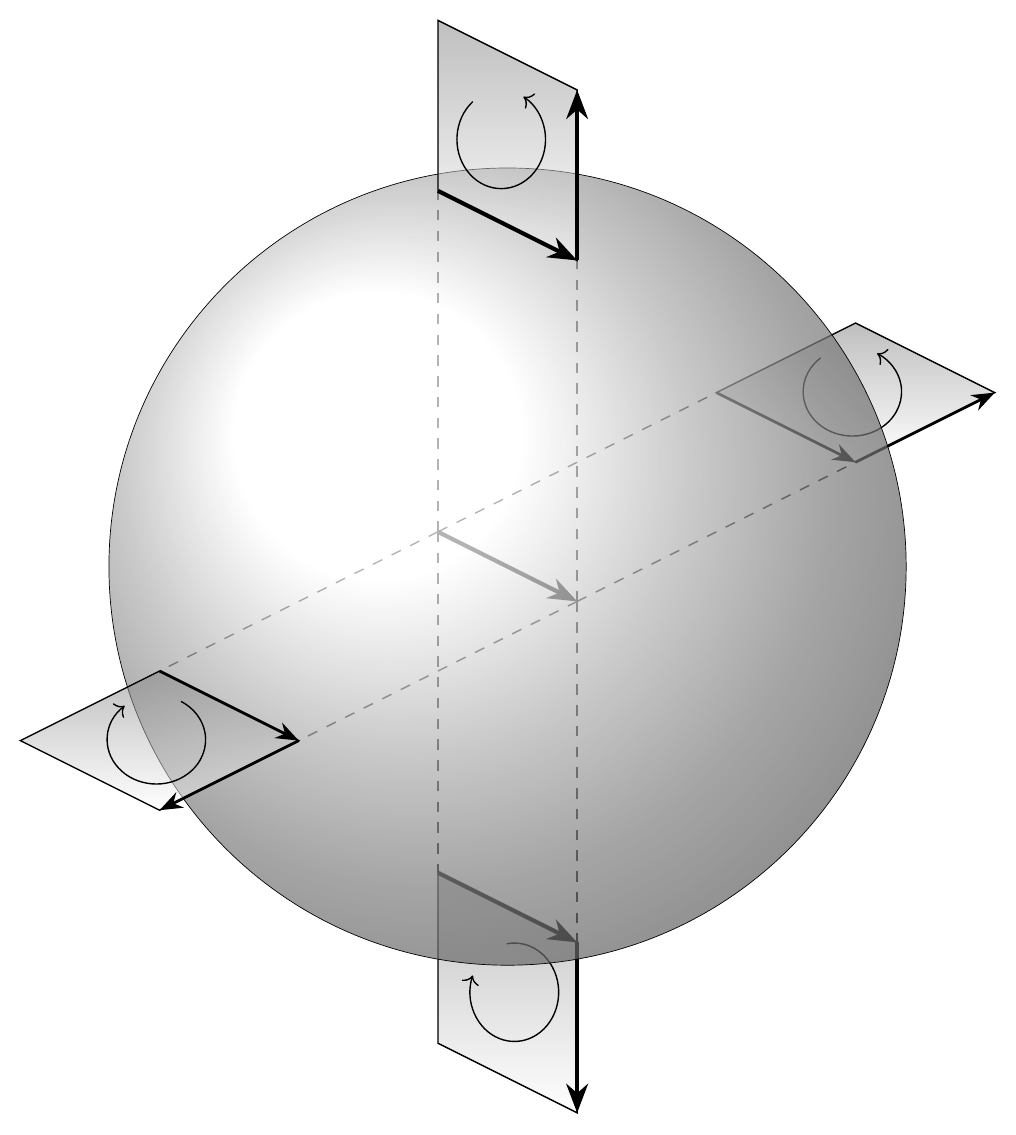}
	\caption{Inner part of the boundary theorem for the volume $M^3$. A bivector field orthogonal to the two-dimensional boundary, equaling to the divergence vector $\boldsymbol{\nabla} \cdot \check{F}$ within.}
	\label{3Dinner}
\end{figure}

For a bivector field $\check{F}$ and volume $M^3=V$ with an oriented trivector volume element $d\tilde{V}=I_3 |d\mathbf{x}^3|$ we can use Eqs.\ \eqref{com} and \eqref{tangvecd} to see that
\begin{equation}
	\boldsymbol{\partial} \check{F} = \lim_{|dV| \to 0} \frac{I_3^{-1}}{|dV|} \oiint_{\partial V} (\text{d}\mathbf{x}^2 \cdot \check{F} + \text{d}\mathbf{x}^2\times \check{F} + \text{d}\mathbf{x}^2\wedge \check{F}).
	\label{bitan}
\end{equation}
The left hand side expands to a vector divergence and a trivector curl, $\boldsymbol{\partial}\check{F}=\boldsymbol{\partial} \cdot \check{F} + \boldsymbol{\partial} \wedge \check{F}$, respectively. The surface integral on the right hand side involves a scalar term, bivector term and a quadvector term, the last of which is zero in $\mathbb{G}^3$. The psuedoscalar division $I_3^{-1}$ in front of the integral, gives the dual to each term, matching the scalar $\text{d}\mathbf{x}^2 \cdot \check{F}$ with the curl, and the bivector $\text{d}\mathbf{x}^2\times \check{F}$ with the divergence. In a similar manner as before, in $\mathbb{R}^3$ we then have the following split of the boundary theorem,
\begin{equation}
	\iiint_{M}\text{d}\mathbf{x}^3 \cdot (\boldsymbol{\nabla} \wedge \mathbf{F})= \oiint_{\partial M} d\mathbf{x}^{2}\cdot \mathbf{F},
	\label{outer2}
\end{equation}
and
\begin{equation}
	\iiint_{M}\text{d}\mathbf{x}^3 \cdot (\boldsymbol{\nabla} \cdot \mathbf{F}   ) = \oiint_{\partial M} d\mathbf{x}^{2}\times \mathbf{F} .
	\label{inner2}
\end{equation}
Fig.\  \ref{3Douter} shows a graphical representation of the outer part in this case, Eq.\ \eqref{outer2}, and suggests an interpretation of the effect of a trivector on a two-dimensional surface enclosing it:  A trivector induces rotation/vorticity at every point on a closed boundary of a volume containing it.
\FloatBarrier
\FloatBarrier
Similarly to the case for vector fields, the inner part for bivector fields in Eq.\ \eqref{inner2}, describes the effect of field components perpendicular to the manifold boundary, Fig.\ \ref{3Dinner}. In the same way that a vector field can diverge from a scalar source, a bivector field diverges from a vector source. Reversing the rotation of the bivector field gives convergence toward a vector sink with a negative sign.

The quadvector term of Eq.\ \eqref{bitan} can be made non-zero by embedding $M^3$ within $\mathbb{R}^4$. In that case, multiplication by $I_3^{-1}$ leaves only the fourth component as $I_3^{-1}I_4=-I_3I_4=-\mathbf{e}_{123}\cdot \mathbf{e}_{1234}= \mathbf{e}_4$. This vector is however perpendicular to $M^3$ and corresponds to the term 
\begin{equation}
	\iiint_{M}\text{d}\mathbf{x}^3 \wedge (\boldsymbol{\nabla} \cdot \check{F}   ) = \oiint_{\partial M} d\mathbf{x}^{2}\wedge \check{F} .
	\label{inner2inR4}
\end{equation}
 The structure can be projected into three-dimensions but requires the manifold $M$ to be embedded in a higher-dimensional one. It's worth noting that this property along with the geometrical shape and genus is shared with the Clifford torus \cite{CT_Taimanov,CT_Kilian,CT_Brendle2013}.  
\section{Dualities}
The dual of an $l$-vector in $\mathbb{G}_n$ is an $m$-vector such that $l+m=n$. With this in mind the dualities of the considered cases of the boundary theorem can be explored. In one dimension, scalars and vectors are dual so the dual of  Fig.\ \ref{fund} is a divergence of sorts, where a scalar within the manifold equals vectors orthogonal to the zero-dimensional boundary. In two dimensions, the outer part in Fig.\ \ref{2Dgreen} and inner part in Fig.\ \ref{2Ddiv} are dual, highlighting the duality of rotation and divergence of vector fields for two dimensional spaces.  

In three dimensions, the outer part shown in Fig.\ \ref{3Douter} is dual to the three-dimensional divergence theorem, Eq.\ \eqref{div}. 
The trivector $\boldsymbol{\nabla} \wedge F$ within $M$ is dual to the divergence scalar in $\mathbb{G}_3$.  The tangent bivector field on the two-dimensional boundary is dual to a radially divergent vector field, normally considered in the three-dimensional divergence theorem.

The dual of the inner part in three dimensions is a known curl theorem variant,
$$\iiint_V \boldsymbol{\nabla} \times \mathbf{F} \, dV= - \oiint_S \mathbf{F} \  \times \hat{\mathbf{N}} \, dS \ .$$
However, if the manifold $M$ is embedded within a higher-dimensional one, it is more suitable to consider its dual within 4D space. In four dimensions, bivectors are self dual and vectors are dual to trivectors. In that case, the duality is the outer part in three dimensions, Fig.\ \ref{3Douter}. 

\section{Projections}
The boundary theorem cases in dimension $N$, are orthogonal projections of cases in the adjacent lower dimension, $N-1$. Therefore, the fundamental theorem of calculus in Fig.\ \ref{fund} can be seen as an orthogonal projection (shadow) of Green's curl theorem, Fig.\ \ref{2Dgreen}. Continuing upwards in dimensions, Green's curl theorem can be seen as a single side projection of the trivector in Fig.\ \ref{3Douter}. Surprisingly, the two dimensional divergence theorem in Fig. \ref{2Ddiv} can be seen as a “front view" of the torus in Fig.\ \ref{3Dinner}.
This projection relation, by definition of the boundary theorem, holds true in all dimensions even though pictorial representations on two-dimensional paper break down.

In general, elements of $\mathbb{G}_n$ can be decomposed into components within and orthogonal to the vector space, so that orthogonal projections and rejections can be written in terms of the geometric product using the pseudoscalar, see further \citet{Hestenes_Sobczyk}. 

\section{DISCUSSION}
\begin{figure}[h!]
	\centering
	\begin{minipage}{0.25\textwidth}
		\centering
				\textbf{(a)}
		\includegraphics[height=3.8cm]{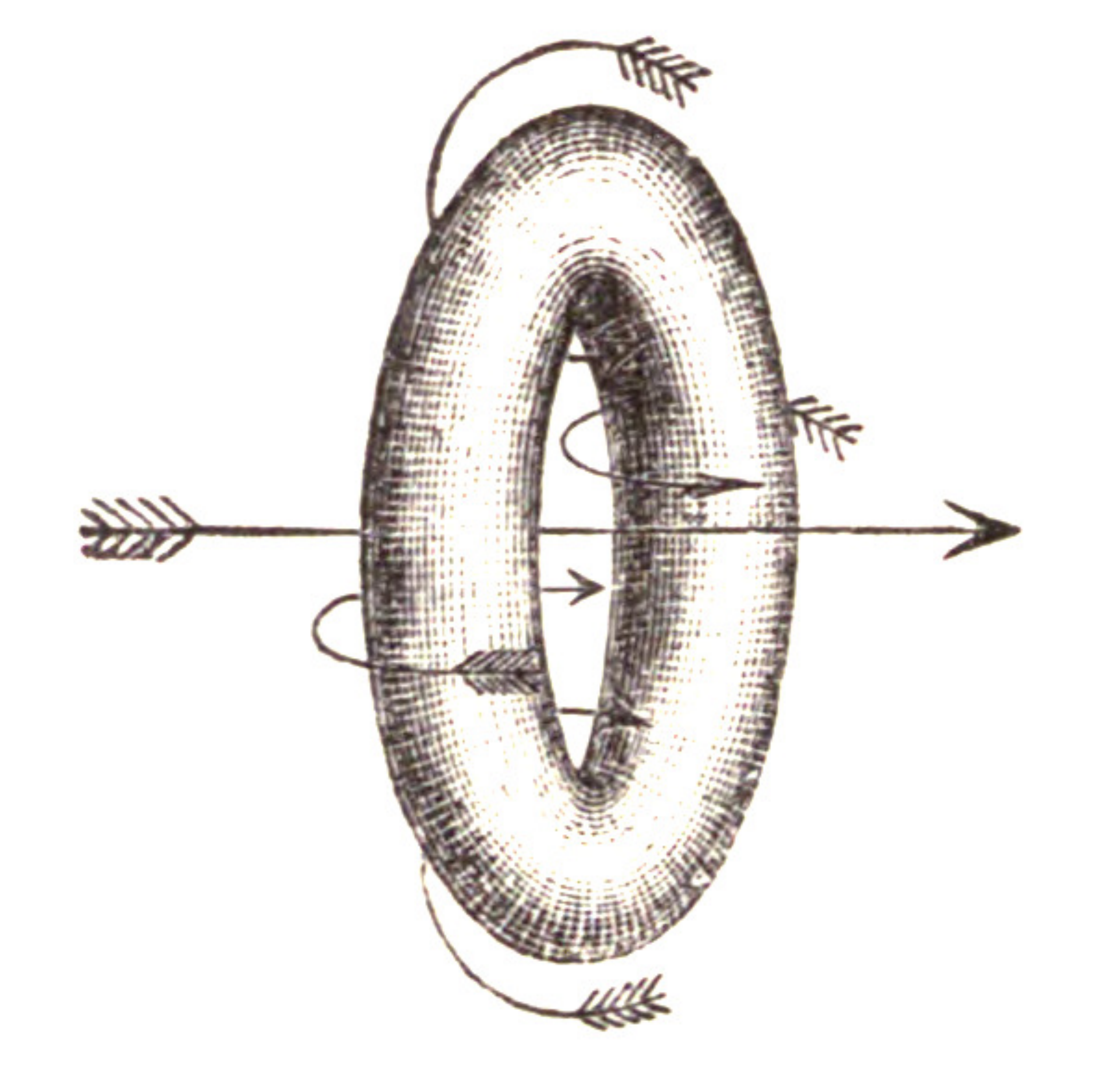}
	\end{minipage}%
	\begin{minipage}{0.25\textwidth}
		\centering
		\textbf{(b)}
		\includegraphics[height=3.8cm]{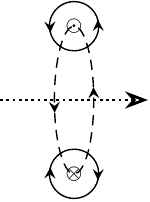}
	\end{minipage}
	\caption{(a) Sketch of a toroidal vortex of a smoke ring \cite{Tait}. (b) Magnetic field (dashed) and magnetic vector potential (solid) around a line current element (dotted). The curl of the magnetic vector potential has the same ring symmetry as the magnetic field, but is shown in a cross-section.  }
	\label{rvs}
\end{figure}

Having uncovered two bivector solutions of the boundary theorem in three-dimensional space, Figs.\ \ref{3Douter} and \ref{3Dinner}, it is natural to consider whether a physical example of these structures exists in nature. Fig.\ \ref{rvs}(a) shows the dynamics of a toroidal vortex ring, as sketched by Tait and published in 1876, in the context of atomic theory \cite{Tait}. Comparison with Fig.\ \ref{3Dinner} shows graphical equivalence. To uncover mathematical equivalence, the vorticity field in fluid dynamics has to be considered in its dual form, as a bivector field   \cite{KOK_MVP}. 
The same dynamic structure can be found if the magnetic vector potential of a current carrying element is considered, from the well-known ring magnetic field induced by it, Fig.\ \ref{rvs}(b).
In both cases, a closed link in the axial field (magnetic/vorticity) of $U(1)$ symmetry is dual to a spin-formation of non-trivial topology in the surrounding medium, hinting at richer internal degrees of freedom \cite{Leurs08}.

Interestingly, another natural phenomenon that corresponds directly to Fig.\ \ref{3Douter} also comes from electromagnetism, though indirectly. Considering that $I\mathbf{B} = \boldsymbol{\nabla} \wedge \mathbf{A}$, one can see that the magnetic vector potential of a magnetic monopole, is essentially the bivector field configuration of Fig.\ \ref{3Douter}. Since $\boldsymbol{\nabla}\wedge \boldsymbol{\nabla} \wedge \mathbf{A}=0$, this indicates that a magnetic monopole should be represented by a trivector, with a bivector potential. It also becomes evident that a magnetic monopole necessitates a radially divergent spin structure. 
 Although we arrive at this conclusion from purely geometric conditions, this can be shown mathematically as well \cite{GAmonopole, STAmonopole,Pezzaglia1994}. This interpretation should apply equally to the monopole configuration of the Berry curvature emerging near degeneracies in parameter space of Hamiltonians \cite{Berry1984,Gradhand_2012} and implies that the first Chern number \cite{Thouless83} is a pseudoscalar trivector.
Note, that if embedded in four dimensions, trivectors will be dual to vectors. Phenomena currently described by scalars in space and vectors in spacetime, may be better understood as trivectors.\\
\pagebreak

Wide applicability of the Figs.\ \ref{3Douter} and \ref{3Dinner} can be expected. The Pauli matrices form an algebra isomorphic to $\mathbb{G}_3$ and the Dirac matrices are a matrix representation of the Dirac algebra $Cl_{1,3}(\mathbb{C})$, both of which have geometric interpretations within space-time algebra\cite{Hestenes_STA} ($Cl_{1,3}(\mathbb{R})$), which is embedded in  $\mathbb{G}_4$. Bott periodicity, key to the periodic table of topological insulators and superconductors \cite{Periodic_Kitaev}, has been shown to have roots in the Octonion division algebra \cite{Octonion_Baez} which has multiple representations in geometric algebra \cite{Octonion_Lasenby}.
Furthermore, every Lie algebra has been shown to have a representation as a bivector algebra \cite{Lie_Spin}. For example, the bivectors of $\mathbb{G}_3$ span the SU(2) Lie group\cite{Hestenes_STA}. This opens up the possibility for visualizing interactions with Lie spin group symmetries such as spin dynamics, electroweak interactions and Berry connections.\\

\section{CONCLUSIONS}
Using geometric algebra, foundational cases of the generalized Stokes' theorem have been visualized and their relations to the standard theorems of vector calculus clarified. From considering bivector valued fields, little explored instances of the theorem emerge which are found to have natural manifestations of non-trivial topology. Dualities and projection relations have been examined. The diagrams presented have wide applicability and for the bivector field cases, may serve in aiding visualization of interactions with spin group symmetry.\\

\begin{acknowledgements}
	This research was supported by the Icelandic Research Fund, grant no. 206568-051. We are grateful to Chris Doran, David Hestenes, Luke Burns and Timo Alho for discussions.
\end{acknowledgements}

\bibliographystyle{apsrev4-1}
\bibliography{heimildir}

\end{document}